%% file: cdc_final_arxiv.tex
\newif\ifconfver
\confverfalse     

\newif\ifplainver  
\plainvertrue

\ifplainver
    \confverfalse   
\fi

\ifconfver

\documentclass[conference]{IEEEtran}

\IEEEoverridecommandlockouts

\else
    \ifplainver
        \documentclass[11pt]{article}
        \usepackage{fullpage}
    \else
        \documentclass[11pt,draftcls,onecolumn]{IEEEtran}
    \fi
\fi

\usepackage{fullpage,graphicx,algorithm,algorithmic,bm,amsmath,amsthm,amssymb,color,shortcuts_OPT,hyperref,cite,enumitem,microtype,url}

\newcounter{algsubstate}


\def\cu#1{{\color{blue}#1}}  

\newtheorem{Fact}{Fact}
\newtheorem{Lemma}{Lemma}
\newtheorem{Prop}{Proposition}
\newtheorem{Theorem}{Theorem}

\newtheorem{Assumption}{H\!\!}

\hypersetup{
  colorlinks   = true, 
  urlcolor     = blue, 
  linkcolor    = blue, 
  citecolor    = blue   
}

\makeatletter
\newcommand{\pushright}[1]{\ifmeasuring@#1\else\omit\hfill$\displaystyle#1$\fi\ignorespaces}
\newcommand{\pushleft}[1]{\ifmeasuring@#1\else\omit$\displaystyle#1$\hfill\fi\ignorespaces}
\makeatother

\usepackage{pgfplots}
\usetikzlibrary{arrows,shapes,calc,tikzmark,backgrounds,matrix,decorations.markings}

\pgfplotsset{compat=1.3}

\usepackage{relsize}
\tikzset{fontscale/.style = {font=\relsize{#1}}
    }

\definecolor{lavander}{cmyk}{0,0.48,0,0}
\definecolor{violet}{cmyk}{0.79,0.88,0,0}
\definecolor{burntorange}{cmyk}{0,0.52,1,0}
\definecolor{asuorange}{rgb}{1,0.699,0.0625}
\definecolor{asured}{rgb}{0.598,0,0.199}
\definecolor{asuborder}{rgb}{0.953,0.484,0}
\definecolor{asugrey}{rgb}{0.309,0.332,0.340}
\definecolor{asublue}{rgb}{0,0.555,0.836}
\definecolor{asugold}{rgb}{1,0.777,0.008}

\tikzstyle{server}=[draw, regular polygon, regular polygon sides=4, black!80, fill=black!40,
                       text=violet,minimum width=10pt]
\tikzstyle{worker}=[draw,circle, asublue!80!white, fill = asublue!50!white,
                       text=violet,minimum width=10pt]

    \makeatletter
    \def\multilimits@{\bgroup
  \Let@
  \restore@math@cr
  \default@tag
 \baselineskip\fontdimen10 \scriptfont\tw@
 \advance\baselineskip\fontdimen12 \scriptfont\tw@
 \lineskip\thr@@\fontdimen8 \scriptfont\thr@@
 \lineskiplimit\lineskip
 \vbox\bgroup\ialign\bgroup\hfil$\m@th\scriptstyle{##}$\hfil\crcr}
    \def\Sb{_\multilimits@}
    \def\endSb{\crcr\egroup\egroup\egroup}
\makeatother

\newtheoremstyle{t}         
    {\baselineskip}{2\topsep}      
    {\rm}                   
    {0pt}{\bfseries}  
    {}                      
    { }                      
    {\thmname{#1}\thmnumber{#2}.}

\theoremstyle{t}

\usepackage{environ}
\NewEnviron{killcontents}{}

\makeatletter
\DeclareRobustCommand*\cal{\@fontswitch\relax\mathcal}
\makeatother

\def\cu#1{{\color{black}#1}} 

\begin{document}
\title{Resilient Distributed Optimization Algorithms for Resource Allocation}
\author{C\'{e}sar A.~Uribe$^\dagger$, Hoi-To Wai$^\dagger$, Mahnoosh Alizadeh\thanks{CAU and HTW have contributed equally. CAU is with LIDS, MIT, Cambridge, MA, USA. HTW is with Dept.~of SEEM, CUHK, Shatin, Hong Kong. MA is with Dept.~of ECE, UCSB, Santa Barbara, CA, USA. This work is partially supported by UCOP Grant LFR-18-548175 and CUHK Direct Grant \#4055113. E-mails: \url{cauribe@mit.edu}, \url{htwai@se.cuhk.edu.hk}, \url{alizadeh@ucsb.edu}}}
\date{\today}

\maketitle

\begin{abstract}
Distributed algorithms provide \cu{flexibility } over centralized algorithms for resource allocation problems, \cu{e.g., cyber-physical systems.} However, the distributed nature of these algorithms often makes the systems susceptible to man-in-the-middle attacks, especially when messages are transmitted between price-taking agents and \cu{a} central coordinator. We propose a resilient strategy for distributed algorithms under the framework of primal-dual distributed optimization. We formulate a robust optimization model that accounts for Byzantine attacks on the communication channels between agents and coordinator. We propose a resilient primal-dual algorithm using state-of-the-art robust statistics methods. The proposed algorithm is shown to converge to a neighborhood of the robust optimization model, where  the neighborhood's radius is proportional to the fraction of attacked channels. 
\end{abstract} 

\section{Introduction}

Consider the following multi-agent optimization problem involving 
the average of parameters in the constraints:
\beq \label{eq:const}
\begin{array}{rl}
\ds \min_{ \prm_i \in \RR^d, \forall i } & U(\prm) \eqdef \frac{1}{N} \sum_{i=1}^N U_i ( \prm_i )\\
{\rm s.t.} & 
g_t \left( \frac{1}{N} \sum_{i=1}^N \prm_i \right) \leq 0,~t=1,...,T,\vspace{.1cm}\\
& \prm_i \in \Cset_i,~i=1,...,N,
\end{array}
\eeq
where both $U_i : \RR^d \rightarrow \RR$ and $g_t : \RR^d \rightarrow \RR$ are
 continuously differentiable, convex functions, and
$\Cset_i$ is a compact convex set in $\RR^d$. We let ${\bm 0} \in \Cset_i$ and 
\beq \label{eq:Rbd}
\max_{ \prm, \prm' \in \Cset_i } \| \prm - \prm' \| \leq R,~i=1,...,N,
\eeq
such that $R$ is an upper bound on the diameters of $\Cset_i$.

Problem~\eqref{eq:const} arises in many resource allocation problems
with a set of potentially \emph{nonlinear} constraints on the amount of allowable resources,
see Section~\ref{sec:ex} for a detailed exploration. 


We consider a system where there exists a central coordinator and $N$ agents. 
In this context, the function $U_i( \prm_i )$ and parameter $\prm_i$ are the utility of the $i$th agent and the resource controlled by agent $i$, respectively. As the agents work independently, it is desirable to design algorithms that allow the $N$ agents to  solve \eqref{eq:const} cooperatively through communication with the central coordinator. Among others, the primal-dual optimization \cu{methods} \cite{koshal2011multiuser} have been advocated as \cu{they} naturally give rise to decomposable algorithms that favor distributed implementation \cite{palomar2006tutorial}. 
In addition to their practical success, these methods are supported by strong theoretical guarantees where fast convergence to an optimal solution of \eqref{eq:const} is well established. 
However, the distributed nature of these methods also exposes the system to \cu{vulnerabilities} not faced by traditional centralized systems. Precisely, existing algorithms assume the agents, and the communication links between central server and agents, to be \emph{completely trustworthy}. \cu{However,} an attacker can take over \cu{a} sub-system operated by the agents, \cu{and} deliberately edit the messages in these communication links, \ie a Byzantine attack. This \cu{might} result in an unstable system \cu{with possible} damages to hardware \cu{and the system overall}. 

In this paper, we propose strategies for securing primal-dual distributed algorithm\cu{s}, e.g., in \cite{koshal2011multiuser}, tailored to solving a relaxed version of the resource allocation problem \eqref{eq:const}. 
A key observation is that the \cu{existing} algorithm\cu{s} \cu{depend} on reliably computing the average of \cu{a} set of parameter vectors, $\{ \prm_i \}_{i=1}^N$, transmitted by the agents. As a remedy, we apply robust statistics techniques as a subroutine, therefore proposing a \emph{resilient} distributed algorithm that is proven to converge
 to a neighborhood of the optimal solution of a \cu{robust} version of \eqref{eq:const}.

\cu{Vulnerabilities of} various types of distributed algorithms have been identified and addressed in a number of recent studies. Relevant examples are \cite{sundaram2011distributed,pasqualetti2012consensus,gentz2016data,sundaram2018distributed,chen2018resilient} which study secure decentralized algorithms on a  general network topology but consider consensus-based optimization models. Moreover, \cite{feng2014distributed,yin2018byzantine,alistarh2018byzantine} consider a similar optimization architecture as this paper, yet they focus on securing distributed algorithms for machine learning tasks which assumes i.i.d.~functions, a fundamentally different setting from the current paper.
Our work is also related to the literature on robust statistics \cite{Donoho1983,huber2011robust}, and particularly, with the recently rekindled research efforts on high dimensional robust statistics \cite{minsker2015geometric,diakonikolas2016robust,Steinhardt2017}.
These works will be the working horse for our attack resilient algorithm.

Our contributions and organization are as follows. First, we derive a formal model for attack resilient resource allocation via a conservative approximation for the robust optimization problem [cf.~Section~\ref{sec:prob}]. Second, we apply and derive new robust estimation results to secure  distributed resource allocation algorithms [cf.~Section~\ref{sec:robust_stat}]. Third, we provide a non-asymptotic convergence guarantee of the proposed attack resilient algorithm [cf.~Section~\ref{sec:anal}]. In particular, our algorithm is shown to converge to a ${\cal O}(\alpha^2)$ neighborhood to the optimal solution of \eqref{eq:const}, where $\alpha \in [0,\frac{1}{2})$ is the fraction of attacked links.

\textbf{Notations}. Unless otherwise specified, $\| \cdot \|$ denotes the standard Euclidean norm. For any $N \in \NN$, $[N]$ denotes the finite set $\{1,...,N\}$. 
%

\subsection{Motivating Examples} \label{sec:ex}
Our set-up here can be employed in a wide range of optimization problems for resource allocation and networked control in multi-agent systems,  e.g., in the pioneering example of  congestion control in data networks \cite{kelly,low}; in determining the optimal price of electricity and enabling more efficient demand supply balancing (a.k.a. demand response) in smart power distribution systems \cite{mohsenian,lina}; in managing user transmit powers and data rates in wireless cellular networks \cite{mung}; in determining optimal caching policies by content delivery networks \cite{cache}; in optimizing power consumption in wireless sensor networks  with energy-restricted batteries \cite{wsn1,wsn2}; and in designing congestion control systems in urban traffic networks \cite{traffic}. 
These examples would have different utility functions and constraint sets that can be handled through our general formulation in \eqref{eq:const}. For example, in the power/rate control problem in data networks, the cost functions are usually  logarithmic functions associated with rate $\theta_i$, e.g., $U_i(\theta_i) = -\beta_i \log (\theta_i)$. In demand response applications in power distribution systems, the utilities capture the users' benefits from operating their electric appliances under different settings. For example, we can capture the cost function of temperature $\theta_i$ controlled by a price-responsive air conditioner
as $U_i(\theta_i) = b_i (\theta_i- \theta_{\rm{ comf}})^2 - c_i$ \cite{lina}. In terms of constraints, our general nonlinear constraint formulation can not only capture common linear resource constraints such as link capacity in data networks \cite{kelly,low}, but  can also handle important non-linear constraints arising in many different applications. For example, in radial power distribution systems, nonlinear convexified power flow constraints can be included for distributed demand response optimization (to see a description of distribution system power flow constraints, see, e.g., \cite{lavaei2014geometry,taylor2015convex}). This can enable our algorithm to perform demand supply balancing in power disribution systems in a  {\it distributed}  and {\it resilient} fashion.

\section{Primal-dual Algorithm for Resource Allocation} 
This section reviews the basic primal\cu{-}dual algorithm for resource allocation.
Let $\bm{\lambda} \in \RR_+^T$ be the dual variable. We consider the Lagrangian
function of \eqref{eq:const}:\vspace{-.1cm}
\beq \notag
\begin{split}
& {\cal L} ( \{ \prm_i \}_{i=1}^N; \bm{\lambda} ) \eqdef  
 \frac{1}{N} \sum_{i=1}^N U_i ( \prm_i ) + \sum_{t=1}^T \lambda_t \!~ g_t \Big( \frac{1}{N} \sum_{i=1}^N \prm_i \Big).
\end{split}\vspace{-.1cm}
\eeq
Assuming strong duality holds (e.g., under the Slater's condition), solving problem~\eqref{eq:const} is equivalent to solving its dual
problem:
\beq \label{eq:pd} \tag{P}
\max_{ \bm{\lambda} \in \RR_+^T } ~
\min_{ \prm_i \in \Cset_i, \forall i}~ {\cal L} ( \{ \prm_i \}_{i=1}^N; \bm{\lambda} ) .
\eeq
For a given $\bm{\lambda}$, the inner minimization of \eqref{eq:pd} is 
known as the Lagrangian relaxation of \eqref{eq:const}, which can be interpreted as a  \emph{penalized} resource allocation problem \cite{lina}.

In a distributed setting, the goal is to solve \eqref{eq:const} where the 
agents only observe a \emph{pricing signal} received from the central coordinator,
and this pricing signal is to be updated iteratively at the central coordinator.
As suggested in \cite{koshal2011multiuser}, we apply the primal\cu{-}dual algorithm (PDA) to a regularized version of \eqref{eq:pd}. Let us define
\beq \label{eq:reg}
\begin{split}
& {\cal L}_{ \upsilon } ( \{ \prm_i \}_{i=1}^N; \bm{\lambda} ) \eqdef \\
& \textstyle {\cal L} ( \{ \prm_i \}_{i=1}^N; \bm{\lambda} ) + \frac{\upsilon}{2 N} \sum_{i=1}^N \| \prm_i \|^2 - \frac{\upsilon}{2} \| \bm{\lambda} \|^2,
\end{split}
\eeq
such that ${\cal L}_{ \upsilon } (\cdot)$ is $\upsilon$-strongly convex and
$\upsilon$-strongly concave in $\{ \prm_i \}_{i=1}^N$ and $\bm{\lambda}$, respectively. Let 
$k \in \ZZ_+$ be the iteration index, $\gamma > 0$ be the step sizes, the PDA recursion is described by:
\begin{subequations} \label{eq:pda}
\begin{align}
& \prm_i^{(k+1)} = \label{eq:pd_prm} \\
& ~~~~{\cal P}_{ \Cset_i } \big( \prm_i^{(k)} - \gamma \!~ \grd_{ \prm_i } 
{\cal L}_{ \upsilon } ( \{ \prm_i^{(k)} \}_{i=1}^N; \bm{\lambda}^{(k)} ) \big),\forall~i \in [N] \notag \\[.1cm]
& \bm{\lambda}^{(k+1)} = \big[ \bm{\lambda}^{(k)} + \gamma \!~ \grd_{ \bm{\lambda} } {\cal L}_{ \upsilon } ( \{ \prm_i^{(k)} \}_{i=1}^N; \bm{\lambda}^{(k)} ) \big]_+ 
\label{eq:pd_lam}
\end{align}
\end{subequations}
where ${\cal P}_{\Cset_i}(\cdot)$ is the Euclidean projection operator, $[\cdot]_+$ denotes $\max\{0,\cdot\}$, and the gradients are:
\beq \label{eq:grd_prm}
\begin{split}
& \grd_{ \prm_i } {\cal L}_{ \upsilon } ( \{ \prm_i^{(k)} \}_{i=1}^N; \bm{\lambda}^{(k)}   ) = \textstyle \frac{1}{N} \Big( \grd_{\prm_i} U_i( \prm_i^{(k)} ) + \upsilon \!~ \prm_i^{(k)} \\
& \hspace{1.6cm} \textstyle + \sum_{t=1}^T \lambda_t^{(k)} \grd_{\prm} g_t (\prm) \Big|_{\prm = \frac{1}{N} \sum_{i=1}^N \prm_i^{(k)}}  \Big),\\[-.6cm]
\end{split}
\eeq
\beq \label{eq:grd_lam}
\begin{split}
& \big[ \grd_{ \bm{\lambda} } {\cal L}_{ \upsilon } ( \{ \prm_i^{(k)} \}_{i=1}^N; \bm{\lambda}^{(k)}  ) \big]_t = g_t \Big( {\textstyle \frac{1}{N} \sum_{i=1}^N} \prm_i^{(k)} \Big) - \upsilon\!~\lambda_t^{(k)},
\end{split}
\eeq
for all $i$, $t$.
We denoted $[{\bm x}]_t$ as the $t$th element of ${\bm x} \in \RR^T$. In particular, observe that \eqref{eq:pda} performs a projected gradient descent/ascent on the primal/dual variables. 

From the above, both gradients with respect to (w.r.t.) $\prm_i$ and $\lambda_t$
depend only on the average parameter $\overline{\prm}^{(k)} \eqdef \frac{1}{N} \sum_{i=1}^N \prm_i^{(k)}$. 
We summarize  the primal dual distributed resource allocation (PD-DRA) procedure in Algorithm~\ref{alg:pddra}.
In addition to solving the general problem \eqref{eq:const}, Algorithm~\ref{alg:pddra} also serves as a general solution method to popular resource allocation problems \cite{lina}.

\algsetup{indent=1em}
\begin{algorithm}[tb]
	\caption{PD-DRA Procedure.}\label{alg:pddra}
	\begin{algorithmic}[1]
        \FOR {$k=1,2,...$}
        \STATE \label{step:accum0} \emph{(Message exchanges stage)}:
        \begin{enumerate}[label=(\alph*)]
        \item \label{step:accum} Central coordinator receives $\{ \prm_i^{(k)} \}_{i=1}^N$ from   agents and computes $\overline{\prm}^{(k)}$, $\{ \grd_{\prm} g_t (\overline{\prm}^{(k)}) \}_{t=1}^T$.
        \item Central coordinator broadcasts the vectors $\overline{\prm}^{(k)}$, $\overline{\bm g}^{(k)} \eqdef \sum_{t=1}^T \lambda_t^{(k)} \grd_{\prm} g_t (\overline{\prm}^{(k)}) $ to   agents.
        \end{enumerate}
        \STATE \label{step:cn} \emph{(Computation stage)}:
        \begin{enumerate}[label=(\alph*)]
        \item \label{step:parallel} Agent $i$ computes the update for $\prm_i^{(k+1)}$ according to \eqref{eq:pd_prm} using the received $\overline{\prm}^{(k)}$.
\item The central coordinator computes the update for $\bm{\lambda}^{(k+1)}$ according to \eqref{eq:pd_lam}. 
        \end{enumerate}
        \ENDFOR
	\end{algorithmic} 
\end{algorithm}


As the regularized primal-dual problem is strongly convex/concave in primal/dual variables, 
Algorithm~\ref{alg:pddra} converges linearly to an optimal solution \cite{koshal2011multiuser}. To study
this, let us denote ${\bm z}^{(k)} = ( \{ \prm_i^{(k)} \}_{i=1}^N, \bm{\lambda}^{(k)} )$
as the primal-dual variable at the  $k$th iteration,   \vspace{-.1cm}
\beq
\bm{\Phi}( {\bm z}^{(k)} ) 
\eqdef \left( \begin{array}{c} 
\grd_{\prm} {\cal L}_{\upsilon} ( \{ \prm_i^{(k)} \}_{i=1}^N, \bm{\lambda}^{(k)}) \\
\grd_{\bm{\lambda}} {\cal L}_{\upsilon} ( \{ \prm_i^{(k)} \}_{i=1}^N, \bm{\lambda}^{(k)} ) 
\end{array} \right).\vspace{-.1cm}
\eeq 
\begin{Fact} \cite[Theorem 3.5]{koshal2011multiuser}
Assume that the map $\bm{\Phi}( {\bm z}^{(k)} )$ is $L_{\Phi}$ Lipschitz continuous. 
For all $k \geq 1$, we have
\beq
\| {\bm z}^{(k+1)}  - {\bm z}^\star \|^2 \leq (1 - 2 \gamma \upsilon  + \gamma^2 L_{\Phi}^2 ) \!~ \| {\bm z}^{(k)} - {\bm z}^\star \|^2 \eqs,
\eeq
where ${\bm z}^\star$ is a saddle point to the regularized version of \eqref{eq:pd}. 
\cu{Setting} $\gamma = \upsilon  / L_{\Phi}^2$ gives 
$\| {\bm z}^{(k+1)}  - {\bm z}^\star \|^2 \leq \big( 1 - \upsilon^2 / L_{\Phi}^2 \big) \| {\bm z}^{(k)} - {\bm z}^\star \|^2$, $\forall~k\geq 1$.
\end{Fact}
\vspace{-.3cm}
\section{Problem Formulation} \label{sec:prob}\vspace{-.2cm}
Despite the simplicity and the strong theoretical guarantee, the PD-DRA method is susceptible to attacks on the channels between the central coordinator and the agents, as described below.

\begin{figure}[t]
\centering{\sf 
\ifconfver
\resizebox{.7\linewidth}{!}{\input{./tikz/DRA.tikz.tex}}
\else
\resizebox{.5\linewidth}{!}{\input{./tikz/DRA.tikz.tex}}
\fi}\vspace{-.3cm}
\caption{Illustrating  the PD-DRA algorithm under attack. The uplink for agent $j$ is compromised such that the correct $\prm_j^{(k)}$ is not transmitted to the central node. The up/downlink for agent $i$ are operating properly.}\label{fig:att}\vspace{-.4cm}
\end{figure}
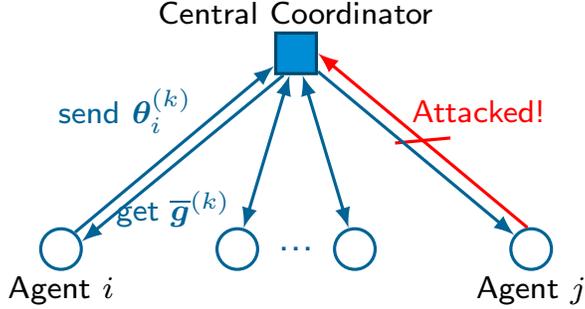

\textbf{Attack Model.}~We consider a situation when \emph{uplink} channels between agents and the central coordinator are compromised [see Fig.~\ref{fig:att}]. Let ${\cal A} \subset [N]$ be the set of \emph{compromised uplink channels}, whose identities are unknown to the central coordinator. We define ${\cal H} \eqdef [N] \setminus {\cal A}$ as the set of trustworthy channels. At iteration $k$, instead of receiving $\prm_i^{(k)}$ from each agent $i \in [N]$ [cf.~Step~\ref{step:accum0}\ref{step:accum}], the central coordinator receives the following messages:\vspace{-.2cm}
\beq \label{eq:signals}
{\bm r}_i^{(k)} = \begin{cases}
\prm_i^{(k)}, & \text{if}~i \in {\cal H},\\
{\bm b}_i^{(k)}, & \text{if}~i \in {\cal A}. 
\end{cases}\vspace{-.2cm}
\eeq
We focus on a Byzantine attack scenario
such that the messages, ${\bm b}_i^{(k)}$,
communicated on the attacked channels can be arbitrary. Under
such scenario, if the central coordinator forms the naive average $\widehat{\prm}^{(k)} = 1/N \sum_{i=1}^N {\bm r}_i^{(k)}$ and computes the gradients 
$\grd g_t( \widehat{\prm}^{(k)} )$ accordingly, this may result in uncontrollable 
error since the deviation $\widehat{\prm}^{(k)} - (1/N) \sum_{i=1}^N \prm_i^{(k)}$ 
can be arbitrarily large. 
It is anticipated that  the PD-DRA method would not provide a solution to the regularized version of \eqref{eq:pd}.



\textbf{Robust Optimization Model.} In light of the Byzantine attack, it is impossible to optimize the original 
problem~\eqref{eq:pd} 
since the contribution from $U_i( \cdot ): i \in {\cal A}$ becomes
unknown to the central coordinator. As a compromise, we focus on optimizing the cost function of agents with trustworthy uplinks and the following robust optimization problem as our target model:
\begin{subequations} \label{eq:robust}
\begin{align}
 \min_{ \begin{subarray}{c} \prm_i \in \Cset_i, i \in {\cal H} \end{subarray}} &~{\textstyle \frac{1}{|{\cal H}|} \sum_{i \in {\cal H}} U_i ( \prm_i )} \\
~\hspace{-.2cm} {\rm s.t.} & \ds
\max_{ \prm_j \in {\cal C}_j, j \in {\cal A}} ~g_t \Big( {\textstyle \frac{1}{N} \sum_{i=1}^N \prm_i} \Big) \leq 0,~\forall~t, \vspace{.1cm} \label{eq:robust_b}
\end{align}
\end{subequations}
note that $\{ \prm_j \}_{j \in {\cal A}}$ is taken away from the decision variables and
we have included \eqref{eq:robust_b} to account for the \emph{worst case} scenario 
for the resource usage of the agents with compromised uplinks.  
This is to ensure that the physical operation limit of the system will not be violated under attack.
Consider the following assumption which will be assumed throughout the paper:
\begin{Assumption} \label{ass:bdd}
For all $\prm \in \RR^d$, the gradient of $g_t$ is bounded with $\| \grd g_t( \prm ) \| \leq B$ and is $L$-Lipschitz continuous.
\end{Assumption}
We define\vspace{-.2cm}
\beq
\overline{g}_t ( \prm ) \eqdef g_t( \prm ) + {\textstyle \frac{|{\cal A}|}{N}} \big( R B + {\textstyle \frac{1}{2}} L R^2 \big),
\eeq 
\begin{Lemma} \label{lem:conservative}
Under H\ref{ass:bdd}. 
The following problem yields a \emph{conservative} approximation of \eqref{eq:robust}, \ie its feasible set is 
a subset of the feasible set of \eqref{eq:robust}:
\beq \label{eq:reformulate_s}
\begin{array}{rl}
\ds \min_{ \prm_i \in \Cset_i,  i \in {\cal H}} & \frac{1}{|{\cal H}|} \sum_{i \in {\cal H}} U_i ( \prm_i ) \vspace{.1cm} \\
{\rm s.t.} & \ds \overline{g}_t \left( {\textstyle \frac{1}{N} \sum_{i \in {\cal H}} \prm_i } \right) \leq 0,~\forall~t \in [T],
\end{array}
\eeq
\end{Lemma}
Similar to PD-DRA, we define the regularized Lagrangian function of \eqref{eq:reformulate_s} as:
\beq
\begin{split}
& \overline{\cal L}_{\upsilon} ( \{ \prm_i \}_{i \in {\cal H}}; \bm{\lambda} ; {\cal H} )  \\
& \eqdef {\textstyle \frac{1}{|{\cal H}|} \sum_{i \in {\cal H} }} U_i ( \prm_i ) + {\textstyle \sum_{t=1}^T} \lambda_t \!~ \overline{g}_t \left( {\textstyle \frac{1}{N} \sum_{i \in {\cal H}} \prm_i} \right) \\
& \hspace{.5cm} \textstyle + \frac{\upsilon}{2 |{\cal H}|} \sum_{i \in {\cal H}} \| \prm_i \|^2 - \frac{\upsilon}{2} \| \bm{\lambda} \|^2.
\end{split}
\eeq
Again, the regularized Lagrangian function is $\upsilon$-strongly convex and concave in
$\prm$ and $\bm{\lambda}$, respectively. 

Our main task is to tackle the following modified problem of \eqref{eq:pd} under Byzantine attack on (some of) the uplinks:
\beq
\label{eq:pdp} \tag{P'}
\max_{ \bm{\lambda} \in \RR_+^T } 
\min_{ \prm_i \in  \Cset_i, \forall i \in {\cal H}}~ \overline{\cal L}_{\upsilon} ( \{ \prm_i \}_{i \in {\cal H}}; \bm{\lambda}; {\cal H} ) ,
\eeq
and we let $\widehat{\bm z}^\star = (\widehat{\prm}^\star, \widehat{\bm{\lambda}}^\star)$ be the optimal
solution to \eqref{eq:pdp}.
Notice that \eqref{eq:pdp} bears a similar form as \eqref{eq:pd} and thus  one may  apply the PD-DRA method to the former naturally. However, such application requires the central coordinator to compute the sample average
\beq \label{eq:truemean}
\textstyle
\overline{\prm}^{(k)}_{\cal H} \eqdef \frac{1}{|{\cal H}|} \sum_{ i \in {\cal H} } \prm_i^{(k)},
\eeq 
at each iteration.
However, the above might not be computationally feasible under the  attack model,
since the central coordinator is oblivious to the identity of ${\cal H}$. 
This is the main objective in the design of our scheme.


\section{Robust Distributed Resource Allocation} \label{sec:robust_stat}
In this section, we describe two estimators for approximating $\overline{\prm}^{(k)}_{\cal H}$ [cf.~\eqref{eq:truemean}] from the received messages \eqref{eq:signals}  without knowing the identity of links in ${\cal H}$. 
To simplify notations, we define $\alpha \geq |{\cal A}| / N$ as a known upper bound to the fraction of compromised channels  and assume $\alpha < 1/2$ where less than half of the channels are compromised.

As discussed after \eqref{eq:truemean},
the problem at hand is \emph{robust mean estimation}, whose applications to robust distributed optimization has been considered in the machine learning literature \cite{yin2018byzantine,alistarh2018byzantine,diakonikolas2016robust} under the assumption that the `trustworthy' signals are drawn i.i.d.~from a Gaussian distribution. Our setting is different since the signals $\prm_i^{(k)}$, $i \in {\cal H}$ are variables from the previous iteration whose distribution is non-Gaussian in general. 
Our analysis will be developed without such assumption on the distribution.   

We first consider a simple median-based estimator applied to each coordinate $j=1,...,d$. First, define the coordinate-wise median as:
\beq 
\big[ {\prm}_{\sf med}^{(k)} \big]_j = {\sf med}\big( \{ [{\bm r}_i^{(k)} ]_j \}_{i=1}^N \big),
\eeq
where ${\sf med}(\cdot)$ computes the median of the operand. Then, our estimator is computed as the mean of the nearest $(1-\alpha)N$ neighbors of $\big[ {\prm}_{\sf med}^{(k)} \big]_j$. To formally describe this, let us define:
\beq
{\cal N}_j^{(k)} = \{ i \in [N] : \big| \big[ {\bm r}_i^{(k)} -  {\prm}_{\sf med}^{(k)} \big]_j \big| \leq r_j^{(k)} \},
\eeq where
$r_j^{(k)}$ is chosen as $|{\cal N}_j^{(k)}| = (1-\alpha)N$. 
Our estimator is:
\beq \label{eq:median} \textstyle
[ \widehat{\prm}^{(k)}_{{\cal H}} ]_j = \frac{1}{(1-\alpha)N} \sum_{i \in {\cal N}_j^{(k)}}  [{\bm r}_i^{(k)} ]_j. 
\eeq
The following bounds the performance of \eqref{eq:median}.
\begin{Prop} \label{lem:median}
Suppose that 
$\max_{i \in {\cal H}} \big\| \prm_i^{(k)} - \overline{\prm}_{\cal H}^{(k)} \big\|_\infty \leq r$,
then for any $\alpha \in (0, \frac{1}{2})$, it holds that 
\beq \label{eq:median_bd}
\big\| \widehat{\prm}^{(k)}_{\cal H} - \overline{\prm}_{\cal H}^{(k)} \big\| \leq \frac{\alpha}{1-\alpha} \Big( 2 + \sqrt{\frac{(1-\alpha)^2}{1-2\alpha}} \Big) r \sqrt{d} \eqs.\vspace{-.1cm}
\eeq
\end{Prop}
Under mild assumptions, the condition $\max_{i \in {\cal H}} \big\| \prm_i^{(k)} - \overline{\prm}_{\cal H}^{(k)} \big\|_\infty \leq r$ can be satisfied with $r = \Theta(R)$, as implied by the compactness of $\Cset_i$ [cf.~\eqref{eq:Rbd}].  
Moreover, for sufficiently small $\alpha$, the right hand side on \eqref{eq:median_bd} can be approximated by ${\cal O}(\alpha R \sqrt{d})$.
However, this median-based estimator may perform poorly for large $\alpha$ (especially when $\alpha \rightarrow \cu{{1}/{2}}$) or dimension $d$. For these situations, a more sophisticated estimator is required, as detailed next.


\begin{algorithm}[t]
	\caption{Recovering the mean of a set \cite{Steinhardt2017}.}
	\begin{algorithmic}[1]\label{algo:1}
		\renewcommand{\algorithmicrequire}{\textbf{Initialize:}}
		\renewcommand{\algorithmicensure}{\textbf{Output:}}
		\STATE \textbf{Input}: $\alpha$, $\prm_i^{(k)}$, $c_i=1$ for all $i=1,\hdots,N$, and $\mathcal{B} = \{1,\hdots,N \}.$\\
		\STATE Set $X_{\mathcal{B}} = [\cdots \prm_j^{(k)}  \cdots ]^\top$ for $j \in \mathcal{B}$ as the concatenated data matrix.
		\STATE Let $Y \in \mathbb{R}^{d\times d}$ and $W \in \mathbb{R}^{\mathcal{B}\times\mathcal{B}}$ be the maximizer/minimizer of the saddle point problem
		{\small
		\begin{align*}
		\max_{\substack{Y \succeq 0, \\ \text{tr}(Y) \leq 1}}
		\min_{\substack{ 0\leq {W_{ij}} ,\\ W_{ij} \leq \frac{4{-}\alpha}{\alpha(2{+}\alpha){\color{red}n}}, \\ \sum_j W_{ji}=1}} \hspace{-.2cm}
		\sum_{i\in\mathcal{B}} c_i (\prm_i^{(k)} {-} X_{\mathcal{B}}w_i)^\top Y(\prm_i^{(k)} {-} X_{\mathcal{B}}w_i)
		\end{align*}\vspace{-.2cm}
	} \label{algo:line_main}
		\STATE Let $\tau_i^{*} = (\prm_i^{(k)} - X_{\mathcal{B}}w_i)^\top Y(\prm_i^{(k)} - X_{\mathcal{B}}w_i)$. 
		\IF {$\sum_{i\in\mathcal{B}} c_i \tau_i^{*}>4n\sigma^2$} \label{algo:lineopt}
		\STATE For $i\in \mathcal{B}$, replace $c_i$ with $\big(1-\frac{\tau^{*}_i}{ \max_{j \in\mathcal{B}}\tau^{*}_j } \big)c_i$.
		\STATE For all $i$ with $c_i < \frac{1}{2}$, remove from $\mathcal{B}$.
		\STATE Go back to Line~\ref{algo:line_main}.
		\ENDIF
		\STATE Set $W_1$ as the result of zeroing out all singular values of $W$ that are greater than $0.9$.
		\STATE Set $Z=X_\mathcal{B} W_0$ where $W_0 = (W {-} W_1)(I{-}W_1)^{\text{-}1}$.
		\IF {$\text{rank}(Z)=1$}
		\STATE \textbf{Output:} $\widehat{\prm}_{\cal H}^{(k)}$ as average of the columns of $X_\mathcal{B}$.
		\ELSE
		\STATE \textbf{Output:} $\widehat{\prm}_{\cal H}^{(k)}$ as a column of $Z$ at random.
		\ENDIF
	\end{algorithmic}
 \end{algorithm}

To derive the second estimator,
we apply an auxiliary result from~\cite{Steinhardt2017} which provides an algorithm for estimating $\overline{\prm}_{\cal H}^{(k)}$, as summarized in Algorithm~\ref{algo:1}. We observe:
\begin{Prop}\label{prop:main} \cite[Proposition 16]{Steinhardt2017}
Suppose that $\lambda_{\max}(\frac{1}{|{\cal H}|}\sum_{i\in {\cal H}} (\prm_i^{(k)} - \overline{\prm}_{\cal H}^{(k)})(\prm_i^{(k)}-\overline{\prm}_{\cal H}^{(k)})^\top ) \leq \sigma^2$. For any $\alpha \in[0, \frac{1}{4})$, Algorithm~\ref{algo:1} produces an output such that $\|\overline{\prm}_{\cal H}^{(k)} - \widehat{\prm}_{\cal H}^{(k)}\| = {\cal O} (\sigma \sqrt{\alpha})$.
\end{Prop} 

Again, similar to Proposition~\ref{lem:median}, the required condition above can be satisfied with $\sigma = \Theta(R)$ under mild conditions. Thus,
Proposition~\ref{prop:main} states that Algorithm~\ref{algo:1} recovers 
 $\overline{\prm}_{\cal H}^{(k)}$ up to an error of ${\cal O}( \sqrt{\alpha} R)$. 
Note that this bound is dimension free unlike the median estimator analyzed in Proposition~\ref{lem:median}.
 
 
 
The idea behind Algorithm~\ref{algo:1} is to sequentially identify and remove the subset of points that cannot be re-constructed from the mean of the data points. The solution of the optimization problem in Line~\ref{algo:line_main} measures how well can we recover the data points as an average of the other $|{\cal H}|$ points. The bounded sample variance assumption guarantees that one can re-construct any element in the set ${\cal H}$ from its mean, thus, all such points that introduce a large error, as quantified by $c_i$ can be safely removed. Line~\ref{algo:lineopt} quantifies the magnitude of the optimal point of Line~\ref{algo:line_main}, and if such value is large, such points that introduce a large error are down-weighted.  The process is repeated until the optimal solution of Line~\ref{algo:line_main} is small enough and a low rank approximation of the optimal $W$ can be used to return the sample mean estimate.



\algsetup{indent=1em}
\begin{algorithm}[t]
	\caption{Resilient PD-DRA}\label{alg:robust_pddra}
	\begin{algorithmic}[1]
		\STATE \textbf{Input}: Each agent has initial state $ \prm_i^{(0)}$.
		\FOR {$k=1,2,...$}
		\STATE \label{r_step:accum0} \emph{(At the Central Coordinator)}:
		\begin{enumerate}[label=(\alph*)]
			\item \label{r_step:accum} Receives $\{{\bm r}_i^{(k)}\}_{i=1}^N$, see~\eqref{eq:signals}, from agents.
			\item Computes robust mean  $\widehat{\prm}_{\cal H}^{(k)}$ using the estimator \eqref{eq:median} or Algorithm~\ref{algo:1}.
			\item Broadcasts the vectors $\widehat{\prm}_{\cal H}^{(k)}$ and $\widehat{\bm g}^{(k)}_{\cal H} \eqdef \sum_{t=1}^T \lambda_t^{(k)} \grd_{\prm} \overline{g_t} (\widehat{\prm}_{\cal H}^{(k)}) $ to   agents.
			\item Computes the update for $\bm{\lambda}^{(k+1)}$  with \eqref{eq:dual_pb}. 
		\end{enumerate}
		\STATE \label{r_step:cn} \emph{(At each agent $i$)}:
		\begin{enumerate}[label=(\alph*)]
			\item Agent receives $\widehat{\prm}_{\cal H}^{(k)}$ and $\widehat{\bm g}^{(k)}_{\cal H}$.
			\item \label{r_step:parallel} Agent computes update for $\prm_i^{(k+1)}$ with \eqref{eq:primal_pb}.

		\end{enumerate}
		\ENDFOR
	\end{algorithmic} 
\end{algorithm}

\textbf{Attack Resilient PD-DRA method}.
The above section provides the   enabling tool for developing the  resilient PD-DRA method, which we summarize in Algorithm~\ref{alg:robust_pddra}.
The algorithm behaves similarly as Algorithm~\ref{alg:pddra} applied to \eqref{eq:pdp}, with the exception that the central coordinator is oblivious to ${\cal H}$, and it uses a robust mean estimator to find an approximate average for the signals sent through the trustworthy links. This approximate value is used to compute the new price signals, and sent back to  agents.
In particular, the primal-dual updates are described by 
\beq \label{eq:primal_pb}
\prm_i^{(k+1)} \hspace{-.2cm} = {\cal P}_{\Cset_i} \big( \prm_i^{(k)} \hspace{-.1cm} - \hspace{-.05cm} {\textstyle \frac{\gamma}{N}} \big( 
\widehat{\bm g}^{(k)}_{\cal H} + \grd U_i( \prm_i^{(k)} ) + \upsilon \prm_i^{(k)} \big) \big),
\eeq
\beq \label{eq:dual_pb}
\lambda_t^{(k+1)} = \big[ \lambda_t^{(k)} + \gamma \big( \overline{g_t} ( {\textstyle \frac{|{\cal H}|}{N}} \widehat{\prm}_{\cal H}^{(k)} ) - \upsilon \lambda_t^{(k)} \big) \big]_+. 
\eeq

\begin{Lemma} \label{obs:pda}
Algorithm~\ref{alg:robust_pddra} is  a primal-dual algorithm \cite{koshal2011multiuser} for \eqref{eq:pdp} with perturbed gradients:\begin{subequations} \label{eq:perturbed}
\begin{align} 
\widehat{\bm g}_{ \prm }^{(k)} & = \grd_{\prm} \overline{\cal L}_\upsilon ( \prm^{(k)}; \bm{\lambda}^{(k)} ; {\cal H} ) + {\bm e}_{\prm}^{(k)} , \label{eq:per_prm} \\[.1cm]
\widehat{\bm g}_{ \bm{\lambda} }^{(k)} & =  \grd_{\bm{\lambda}} \overline{\cal L}_\upsilon (  \prm_i^{(k)} ; \bm{\lambda}^{(k)} ; {\cal H} ) + {\bm e}_{\bm{\lambda}}^{(k)},
\label{eq:per_lam}
\end{align}
\end{subequations}
where we have used concatenated variable as $\prm = (\prm_1,...,\prm_N)$ and $\bm{\lambda} = (\lambda_1,...,\lambda_T)$. Under H\ref{ass:bdd} and assuming that $\lambda_t^{(k)} \leq \overline{\lambda}$ for all $k$, we have:
\beq
\|  {\bm e}_{\prm}^{(k)} \| \leq \overline{\lambda} LT \| \widehat{\prm}_{\cal H}^{(k)} - \overline{\prm}_{\cal H}^{(k)} \|, \label{eq:bd_prm}
\eeq
\beq \|{\bm e}_{\bm{\lambda}}^{(k)}\| \leq BT \| \widehat{\prm}_{\cal H}^{(k)} - \overline{\prm}_{\cal H}^{(k)} \|. \label{eq:bd_lambda}\vspace{-.1cm}
\eeq
\end{Lemma}
The assumption $\lambda_t^{(k)} \leq \overline{\lambda}$ can be guaranteed since $\overline{g_t} ( {\textstyle \frac{|{\cal H}|}{N}} \widehat{\prm}_{\cal H}^{(k)} )$ is bounded.\vspace{-.1cm}

\subsection{Convergence Analysis} \label{sec:anal}
%
%
%
Finally, based on Lemma~\ref{obs:pda}, we can analyze the convergence of Algorithm~\ref{alg:robust_pddra}.
Let $\widehat{\bm z}^\star = (\widehat{\prm}^\star, \widehat{\bm{\lambda}}^\star)$
be a saddle point of \eqref{eq:pdp} and define 
\beq
\overline{\bm{\Phi}}( {\bm z}^{(k)} ) 
\eqdef \left( \begin{array}{c} 
\grd_{\prm} \overline{\cal L}_{\upsilon} ( \{ \prm_i^{(k)} \}_{i \in {\cal H}}, \bm{\lambda}^{(k)}; {\cal H} ) \\
- \grd_{\bm{\lambda}} \overline{\cal L}_{\upsilon} ( \{ \prm_i^{(k)} \}_{i \in {\cal H}}, \bm{\lambda}^{(k)} ; {\cal H} ) 
\end{array} \right),
\eeq 
\begin{Theorem} \label{thm:main}
Assume the map $\overline{\bm{\Phi}}( {\bm z}^{(k)} )$ is $L_{\Phi}$-Lipschitz continuous. For Algorithm~\ref{alg:robust_pddra}, for all $k \geq 0$ it holds
\beq \label{eq:rate_E} \begin{split}
& \| {\bm z}^{(k+1)} - \widehat{\bm z}^\star \|^2 \leq \big( 1 - \gamma \upsilon + 2 \gamma^2 L_{\Phi}^2 \big) \| {\bm z}^{(k)} - \widehat{\bm z}^\star \|^2 \\
&+ \big( \frac{4 \gamma}{\upsilon} + 2 \gamma^2  \big) E_k.  \end{split}
\eeq
where $E_k \eqdef \| {\bm e}_{\prm}^{(k)} \|^2 + \| {\bm e}_{\bm{\lambda}}^{(k)} \|^2$ is the total perturbation at iteration $k$. Moreover, if we choose 
$\gamma < {\upsilon} / {2L_{\Phi}^2}$
and $E_k$ is upper bounded by $\overline{E}$ for all $k$, then
\beq \label{eq:upperbd_E}
\limsup_{k \rightarrow \infty} \| {\bm z}^{(k)} - \widehat{\bm z}^\star \|^2 \leq 
\frac{ \frac{4}{\upsilon} + 2 \gamma }{ \upsilon - 2 \gamma L_{\Phi}^2} \overline{E}\cu{.}\vspace{-.2cm}
\eeq
\end{Theorem}
Combining the results from the last subsection,
the theorem shows the desired result that the resilient PD-DRA method converges to
a ${\cal O}(\alpha^2 R^2 d)$ neighborhood of the saddle point of \eqref{eq:pdp}, if the median-based estimator \eqref{eq:median} is used [or ${\cal O}(\alpha R^2)$ if Algorithm~\ref{algo:1} is used], where $\alpha$ is the fraction of attacked uplink channels. Moreover, it shows that the convergence rate to the neighborhood is linear, which is similar to the classical PDA analysis \cite{koshal2011multiuser}.

Interestingly, Theorem~\ref{thm:main} illustrates a trade-off in the choice of the step size $\gamma$
between convergence speed and accuracy.  
In specific, \eqref{eq:rate_E} shows that the rate of convergence 
factor $1 - \gamma \upsilon + 2 \gamma^2 L_\Phi^2$ can be minimized
by setting $\gamma = \upsilon / (4L_\Phi^2)$. 
However, in the meantime, the asymptotic upper bound in \eqref{eq:upperbd_E}
is increasing with $\gamma$ and it can be minimized by setting $\gamma \rightarrow 0$. 
This will be a design criterion to be explored in practical implementations.\vspace{-.1cm}

\section{Conclusions}\vspace{-.1cm}
In this paper, we studied the strategies for securing a primal-dual algorithm for distributed resource allocation. Particularly, we propose a resilient distributed algorithm based on primal-dual optimization and robust statistics. We derive bounds for the performance of the studied algorithm and show that it converges to a neighborhood of a robustified resource allocation problem when the number of attacked channels is small. 

\section*{Acknowledgement} 
The authors would like to thank the anonymous reviewers for feedback, and Mr.~Berkay Turan (UCSB) for pointing out typos in the original submission of this paper.

\iftrue
\ifconfver
\appendices
\else
\appendix
\fi

\section{Proof of Lemma~\ref{lem:conservative}} \label{app:conservative}
Since $g_t$ is $L_t$-smooth, the following holds
\beq \label{eq:some_ineq}
\begin{split}
& \textstyle g_t \big( \frac{1}{N} \sum_{i=1}^N \prm_i \big) 
\leq g_t \big( \frac{1}{N} \sum_{i \in {\cal H}} \prm_i \big) \\
& \textstyle + 
\frac{1}{N} \sum_{j \in {\cal A}} \Big\langle \prm_j, \grd g_t \big( \frac{1}{N} \sum_{i \in {\cal H}} \prm_i \big) \Big\rangle + \frac{L_t}{2 N^2} \big\| \sum_{j \in {\cal A}} \prm_j \big\|^2
\end{split} 
\eeq
Furthermore, observe that the gradient of $g_t$ is uniformly bounded by $B$ and the diameter of ${\cal C}_j$
is $R$, then the right hand side of \eqref{eq:some_ineq} can be upper bounded by
\beq
g_t \big( {\textstyle \frac{1}{N} \sum_{i \in {\cal H}} \prm_i} \big) + {\textstyle \frac{1}{N}} \sum_{j \in {\cal A}} \big( R B + {\textstyle \frac{1}{2}} L R^2 \big)
\eeq
As such, defining $c_t \eqdef \frac{|{\cal A}|}{N}\big( R B + \frac{1}{2} L R^2 \big)$,
it can be seen that
\beq
g_t \big( {\textstyle \frac{1}{N} \sum_{i \in {\cal H}} \prm_i} \big) + c_t \leq 0,~t=1,...,T 
\eeq
implies the desired constraint in \eqref{eq:robust}. 

\section{Proof of Proposition~\ref{lem:median}}\label{app:median}
Fix any $j \in [d]$. The assumption implies that for all $i \in {\cal H}$, one has 
\beq
\big| [ \prm_i^{(k)} - \overline{\prm}_{\cal H}^{(k)} ]_j \big| \leq r.
\eeq
We observe that $|{\cal H}| \geq (1-\alpha) N$. Applying \cite[Lemma 1]{feng2014distributed} shows that the median estimator\footnote{At each coordinate, the median is the geometric median estimator of one dimension in \cite{feng2014distributed}.} satisfies  
\beq
\big| \big[ {\prm}_{\sf med}^{(k)} -  \overline{\prm}_{\cal H}^{(k)} \big]_j \big| \leq (1- \alpha) \sqrt{\frac{1}{1-2\alpha}} r.
\eeq
The above implies that for all $i \in {\cal H}$, we have
\beq
\big| [ \prm_i^{(k)} - {\prm}_{\sf med}^{(k)} ]_j \big| \leq \Big( 1 + \sqrt{\frac{(1-\alpha)^2}{1-2\alpha}} \Big) r.
\eeq
This implies that $r_j^{(k)} \leq \Big( 1 + \sqrt{\frac{(1-\alpha)^2}{1-2\alpha}} \Big) r$ since $|{\cal H}| \geq (1-\alpha)N$.
We then bound the performance of $\widehat{\prm}^{(k)}$:
\beq
\begin{split}
& (1-\alpha)N [ \widehat{\prm}^{(k)} ]_j = \sum_{i \in {\cal N}_j^{(k)}}  [{\bm r}_i^{(k)} ]_j \\
& =
\sum_{i \in {\cal H} }  [{\bm r}_i^{(k)} ]_j - \hspace{-.2cm} \sum_{i \in {\cal H} \setminus {\cal N}_j^{(k)} }  [{\bm r}_i^{(k)} ]_j + \hspace{-.2cm} \sum_{i \in {\cal A} \cap {\cal N}_j^{(k)} }  [{\bm r}_i^{(k)} ]_j,
\end{split}
\eeq
thus
\beq \notag
\begin{split}
& (1-\alpha)N [ \widehat{\prm}^{(k)} - \overline{\prm}_{\cal H}^{(k)} ]_j  \\
& = - \sum_{i \in {\cal H} \setminus {\cal N}_j^{(k)} }  [{\bm r}_i^{(k)} - \overline{\prm}_{\cal H}^{(k)} ]_j + \hspace{-.2cm} \sum_{i \in {\cal A} \cap {\cal N}_j^{(k)} }  [{\bm r}_i^{(k)} - \overline{\prm}_{\cal H}^{(k)} ]_j.
\end{split}
\eeq
Notice that $|{\cal A} \cap {\cal N}_j^{(k)}| \leq \alpha N$ and thus $|{\cal H} \setminus {\cal N}_j^{(k)}| \leq \alpha N$. Gathering terms shows
\beq
| [ \widehat{\prm}^{(k)} - \overline{\prm}_{\cal H}^{(k)} ]_j | \leq \frac{\alpha N}{(1-\alpha)N}
\Big( 2 + \sqrt{\frac{(1-\alpha)^2}{1-2\alpha}} \Big) r.
\eeq
The above holds for all $j \in [d]$. Applying the norm equivalence shows the desired bound.

\section{Proof of Lemma~\ref{obs:pda}} \label{app:obs}
Comparing the equations in \eqref{eq:perturbed} with \eqref{eq:primal_pb}, \eqref{eq:dual_pb}, we identify that
\beq \label{eq:etheta}
\big[ {\bm e}_{\prm}^{(k)} \big]_i = \frac{1}{N}  {\sum_{t=1}^T} \lambda_t^{(k)} \big( \grd_{\prm} \overline{g_t} ( {\textstyle \frac{|{\cal H}|}{N} } \widehat{\prm}_{\cal H}^{(k)} ) -  \grd_{\prm} \overline{g_t} ( {\textstyle \frac{|{\cal H}|}{N} } \overline{\prm}_{\cal H}^{(k)} ) \big) 
\eeq
\beq \label{eq:elambda}
\big[ {\bm e}_{\bm{\lambda}}^{(k)} \big]_t = \overline{g}_t ( {\textstyle \frac{|{\cal H}|}{N} }\widehat{\prm}_{\cal H}^{(k)} ) - \overline{g}_t ( {\textstyle \frac{|{\cal H}|}{N} }\overline{\prm}_{\cal H}^{(k)} ) 
\eeq
where $\big[ {\bm e}_{\prm}^{(k)} \big]_i$ denotes the $i$th block of ${\bm e}_{\prm}^{(k)}$.
Using H\ref{ass:bdd} and the said assumptions, we immediately see that
\beq
\| \big[ {\bm e}_{\prm}^{(k)} \big]_i \| \leq \overline{\lambda} \frac{L T}{N}  \| 
\widehat{\prm}_{\cal H}^{(k)} - \overline{\prm}_{\cal H}^{(k)} \| 
\eeq
which then implies \eqref{eq:bd_prm}.
 H\ref{ass:bdd} implies that $\overline{g}_t$ is $B$-Lipschitz continuous,
therefore 
\beq
| \big[ {\bm e}_{\bm{\lambda}}^{(k)} \big]_t | \leq B \| 
\widehat{\prm}_{\cal H}^{(k)} - \overline{\prm}_{\cal H}^{(k)} \|,
\eeq
which implies \eqref{eq:bd_lambda}.

\section{Proof of Theorem~\ref{thm:main}}
Based on Proposition~\ref{obs:pda}, our idea is to perform a perturbation analysis on the PDA algorithm. 
Without loss of generality, we assume $N=1$ and denote $\prm = \prm_1$. 
To simplify notations, we also drop the subscript, denote the modified and regularized Lagrangian function as ${\cal L} = \overline{\cal L}_{\upsilon}$. 
Furthermore, we denote the saddle point to \eqref{eq:pdp} as ${\bm z}^\star = (\prm^\star, \bm{\lambda}^\star)$. 

Using the fact that $\prm^\star = {\cal P}_{\Cset}( \prm^\star  )= {\cal P}_{\Cset} \big( \prm^\star - \gamma \grd_{\prm} {\cal L} ( \prm^\star, \bm{\lambda}^\star) \big)$, we 
observe that in the primal update:
\begin{align*}
& \| \prm^{(k+1)} - \prm^\star \|^2 \label{eq:ss} \\
& \overset{(a)}{\leq} \| \prm^{(k)} - \prm^\star \|^2 - 2 \gamma 
\langle \widehat{\bm g}_{\prm}^{(k)} - \grd_{\prm} {\cal L} ( \prm^\star, \bm{\lambda}^\star), \prm^{(k)} - \prm^\star \rangle \\
& \hspace{.8cm} + \gamma^2 
\| \widehat{\bm g}_{\prm}^{(k)} - \grd_{\prm} {\cal L} ( \prm^\star, \bm{\lambda}^\star) \|^2 
\end{align*}
where (a)  is due to the projection inequality $\| {\cal P}_{\Cset} ( {\bm x} - {\bm y} ) \| \leq 
\| {\bm x} - {\bm y} \|$. Furthermore, using the Young's inequality, for any $c_0, c_1 > 0$, we have
\begin{align*}
& \| \prm^{(k+1)} - \prm^\star \|^2 \\
& \leq \| \prm^{(k)} - \prm^\star \|^2 \\
& \hspace{.5cm} - 2 \gamma \langle \grd_{\prm} {\cal L}( \prm^{(k)}, \bm{\lambda}^{(k)} ) - \grd_{\prm} {\cal L} ( \prm^\star, \bm{\lambda}^\star), \prm^{(k)} - \prm^\star \rangle \\
& \hspace{.5cm} { + \gamma^2 ( 1 + c_0 ) \|  \grd_{\prm} {\cal L}( \prm^{(k)}, \bm{\lambda}^{(k)} ) - \grd_{\prm} {\cal L} ( \prm^\star, \bm{\lambda}^\star) \|^2} \\
& \hspace{.5cm} {- 2 \gamma \langle {\bm e}_{\prm}^{(k)},  \prm^{(k)} - \prm^\star \rangle + \gamma^2 \big( 1 + \frac{1}{c_0 } \big) \| {\bm e}_{\prm}^{(k)} \|^2} \\
& \leq (1 + 2 c_1 \gamma ) \!~  \| \prm^{(k)} - \prm^\star \|^2 \\
& \hspace{.5cm} - 2 \gamma \langle \grd_{\prm} {\cal L}( \prm^{(k)}, \bm{\lambda}^{(k)} ) - \grd_{\prm} {\cal L} ( \prm^\star, \bm{\lambda}^\star), \prm^{(k)} - \prm^\star \rangle \\
& \hspace{.5cm} + \gamma^2 ( 1 + c_0 ) \|  \grd_{\prm} {\cal L}( \prm^{(k)}, \bm{\lambda}^{(k)} ) - \grd_{\prm} {\cal L} ( \prm^\star, \bm{\lambda}^\star) \|^2  \\
& \hspace{.5cm} + \Big( \frac{2 \gamma}{c_1} + \gamma^2 + \frac{\gamma^2}{c_0 } \Big) \| {\bm e}_{\prm}^{(k)} \|^2 .
\end{align*}
Similarly, in the dual update we get, 
\begin{align*}
& \| \bm{\lambda}^{(k+1)} - \bm{\lambda}^\star \|^2 \\
& \leq \| \bm{\lambda}^{(k)} - \bm{\lambda}^\star \|^2 + \gamma^2 
\| \widehat{\bm g}_{\bm{\lambda}}^{(k)} - \grd_{\bm{\lambda}} {\cal L} ( \prm^\star, \bm{\lambda}^\star) \|^2\\
& \hspace{.5cm} + 2 \gamma 
\langle \widehat{\bm g}_{\bm{\lambda}}^{(k)} - \grd_{\bm{\lambda}} {\cal L} ( \prm^\star, \bm{\lambda}^\star), \bm{\lambda}^{(k)} - \bm{\lambda}^\star \rangle \\
& \leq (1 + 2 c_1 \gamma ) \!~\| \bm{\lambda}^{(k)} - \bm{\lambda}^\star \|^2 \\
& \hspace{.5cm} + 2 \gamma \langle \grd_{\bm{\lambda}} {\cal L}( \prm^{(k)}, \bm{\lambda}^{(k)} ) - \grd_{\bm{\lambda}} {\cal L} ( \prm^\star, \bm{\lambda}^\star), \bm{\lambda}^{(k)} - \bm{\lambda}^\star \rangle \\
& \hspace{.5cm} + \gamma^2 ( 1 + c_0 ) \| \grd_{\bm{\lambda}} {\cal L}( \prm^{(k)}, \bm{\lambda}^{(k)} ) - \grd_{\bm{\lambda}} {\cal L} ( \prm^\star, \bm{\lambda}^\star) \|^2  \\
& \hspace{.5cm}+ \Big( \frac{2 \gamma}{c_1} + \gamma^2  + \frac{\gamma^2}{c_0 } \Big) \| {\bm e}_{\bm{\lambda}}^{(k)} \|^2.
\end{align*}
Summing up the two inequalities gives:
\beq \notag
\begin{split}
& \| {\bm z}^{(k+1)} - {\bm z}^\star \|^2 \\
& \leq (1 + 2 c_1 \gamma ) \!~\| {\bm z}^{(k)} - {\bm z}^\star \|^2 + \Big( \frac{2 \gamma}{c_1} +\gamma^2 +  \frac{\gamma^2}{c_0 } \Big) E_k  \\
& \hspace{.25cm} - 2 \gamma \langle \bm{\Phi}( {\bm z}^{(k)} ) - \bm{\Phi}( {\bm z}^\star ),  {\bm z}^{(k)} - {\bm z}^\star \rangle \\
& \hspace{.5cm} + \gamma^2 ( 1 + c_0 ) \| \bm{\Phi}( {\bm z}^{(k)} ) - \bm{\Phi}( {\bm z}^\star ) \|^2  \\
& \overset{(a)}{\leq}  \Big( 1 + 2 \gamma (c_1  - \upsilon ) + \gamma^2 ( 1+ c_0) L_{\Phi}^2 \Big) \| {\bm z}^{(k)} - {\bm z}^\star \|^2 \\
& \hspace{.5cm} + \Big( \frac{2 \gamma}{c_1} +\gamma^2 + \frac{\gamma^2}{c_0 } \Big) E_k,
\end{split}
\eeq
where (a) uses the strong monotonicity and smoothness of the map $\bm{\Phi}$. 
Setting $c_1 = \upsilon /2$ yields
\beq
\begin{split}
& \| {\bm z}^{(k+1)} - {\bm z}^\star \|^2 \\
& \leq \Big( 1 - \gamma \upsilon + \gamma^2 (1+c_0) L_{\Phi}^2 \Big) \| {\bm z}^{(k)} - {\bm z}^\star \|^2 \\
& \hspace{.5cm} +
\Big( \frac{4 \gamma}{\upsilon} +\gamma^2 + \frac{\gamma^2}{c_0 } \Big) E_k.
\end{split}
\eeq
Observe that we can choose $\gamma$ such that $1 - \gamma \upsilon + \gamma^2 (1+c_0) L_{\Phi}^2 < 1$.
Moreover, the above inequality implies that $\| {\bm z}^{(k)} - {\bm z}^\star \|^2$ evaluates
to
\beq \notag
\begin{split}
& \| {\bm z}^{(k+1)} - {\bm z}^\star \|^2 \\
& \leq (1 - \gamma \upsilon + \gamma^2 (1+c_0) L_{\Phi}^2 )^k \| {\bm z}^{(0)} - {\bm z}^\star \|^2 + \\
& \sum_{\ell=1}^k (1 - \gamma \upsilon + \gamma^2 (1+c_0) L_{\Phi}^2 )^{k-\ell} \Big( \frac{4 \gamma}{\upsilon} + \gamma^2+ \frac{\gamma^2}{c_0 } \Big) E_\ell 
\end{split}
\eeq
If $E_k \leq \overline{E}$ for all $k$, then ${\bm z}^{(k)}$ converges to a neighborhood of ${\bm z}^\star$ of radius 
\beq
\limsup_{k \rightarrow \infty} \| {\bm z}^{(k)} - {\bm z}^\star \|^2 \leq 
\frac{ \frac{4 \gamma}{\upsilon} + \gamma^2 +\frac{\gamma^2}{c_0 } }{\gamma \upsilon - \gamma^2 ( 1+ c_0 ) L_{\Phi}^2} \overline{E}
\eeq
Setting $c_0 =1$ concludes the proof.

\fi

\bibliographystyle{ieeetr}
\bibliography{robust_bib,num_bib,opt_bib}


\end{document}

%% file: tikz/DRA.tikz.tex
\begin{tikzpicture}[scale=0.7,thick]
\tikzset{
  strike through/.style={
    postaction=decorate,
    decoration={
      markings,
      mark=at position 0.5 with {
        \draw[-] (-5pt,-5pt) -- (5pt, 5pt);
      }
    }
  }
}
  \node[server, asublue!50!black, fill=asublue ] (ser) at (0,0) {};
  \foreach \place/\name in {{(-3,-2.5)/a}, {(-.75,-2.5)/b}, {(.75,-2.5)/c}, {(3,-2.5)/d}}
    \node[worker, asublue!70!black, fill=white, opacity=1] (\name) at \place {};
\node[asublue!70!black] at (0,-2.5) {...};
\node [below, yshift=.6cm] at (ser) {{\scriptsize Central Coordinator}};
\node [below, yshift=-.1cm] at (a) {{\scriptsize Agent $i$}};
\node [below, yshift=-.1cm] at (d) {{\scriptsize Agent $j$}};
\draw[>=latex, thick, asublue!70!black,<-] (a.20) to node[auto, yshift=-.75cm,xshift=0.575cm] {\scriptsize get $\overline{\bm g}^{(k)}$} (ser.240); 
\draw[>=latex, thick, asublue!70!black,->] (a.50) to node[auto, xshift=.3cm] {\scriptsize send $\prm_i^{(k)}$} (ser.210); 
\draw[>=latex, thick, asublue!70!black,<->] (b) to (ser); 
\draw[>=latex, thick, asublue!70!black,<->] (c) to (ser); 
\draw[>=latex, thick, red,->, strike through] (d.100) -- node[auto, yshift=.5cm, xshift=1.25cm] {\scriptsize Attacked!} 
(ser.0); 
\draw[>=latex, thick, asublue!70!black,<-] (d) to (ser);  
\end{tikzpicture}